\documentclass{article}

\usepackage{default}
\usepackage{mathrsfs}
\usepackage{amsfonts}
\usepackage{amsmath}
\usepackage{mathbbol}
\usepackage[sans]{dsfont}
\usepackage[mathscr]{euscript}
\usepackage[nospace]{overcite}

\renewcommand{\d}{\mathrm{d}}

\newcommand{\1}{\mathds{1}}

\renewcommand{\H}{\mathscr{H}}

\newcommand{\ds}{\displaystyle}
\newcommand{\ot}{\otimes}
\newcommand{\la}{\left\langle}
\newcommand{\ra}{\right\rangle}

\renewcommand{\i}{\mathrm{i}}
\newcommand{\Ce}{C_\epsilon}
\newcommand{\ad}{a^{\dag}}
\newcommand{\bd}{b^{\dag}}
\newcommand{\lii}{\ell^2(\mathbb{N})}

\newcommand{\mr}{\mathrm}
\newcommand{\norm}[1]{\left\Vert#1\right\Vert}

 \newtheorem{theorem}{Theorem}
 \newtheorem{corollary}{Corollary}
 \newtheorem{lemma}{Lemma}
 \newtheorem{proposition}{Proposition}
 \newtheorem{definition}{Definition}

\title{Quantum Stochastic Dynamics in\\Multi-Photon Optics}
\author{Ricardo Castro Santis \\ Departamento de Matem\'atica Universidad del Bio-Bio and\\ANESTOC, Pontificia Universidad Cat\'olica de Chile}

\begin{document}
\date{}
\maketitle
\begin{abstract}
Multi-photon models are theoretically and experimentally important because in them quantum properly phenomena are verified; as well as squeezed light and quantum entanglement also plays a relevant role in quantum information and quantum communication (see [\citen{ZolG97, WuK87, Jeff-Ou-2006 }]).

In this paper we study a generic model of a multi-photon system with an arbitrary number of pumping and subharmonics fields. This model includes measurement on the system, as could be direct or homodyne detection and we demonstrate the existence of dynamics in the context of \emph{Continuous Measurement Theory of Open Quantum Systems} (see [\citen{Acc78, Bar86LNP,  Bar86PR, BAR-1990-QO, BAR-OQS, BarGSpringer, BarL85JMP, Carm08, Belavkin-99, RCS-LAP, Castro-Barchielli}]) using \emph{Quantum Stochastic Differential Equations} with unbounded coefficients (see [\citen{RCS-LAP, Castro-Barchielli, Fagnola99, Fagnola2, Fagnola3, Fagnola-Wills}]).

\end{abstract}


\section{Introduction}
\subsection{The Creator, Annihilation and Number Operators}

Consider the Hilbert space $\lii$ with its canonical base $\left\{e_k\right\}_{k\in\mathbb{N}}$. The creator and annihilation operators are defined in $\lii$ by:

\begin{equation*}\label{eq:cre&annihi-operators}
\begin{array}{c}
  \mathrm{Dom}(\ad)=\mathrm{Dom}(a)=\left\{\{u_k\}\in\lii;\quad \sum_{k\ge1}k|u_k|\le\infty\right\}\\
\\
   \ad e_k=\sqrt{k+1}\ e_{k+1},\qquad ae_k=\sqrt{k}\ e_{k-1}\qquad ae_0=0
\end{array}
\end{equation*}

The number operator is defined  by:
 \begin{equation*}\label{eq:number-operators}
  \mathrm{Dom}(\ad a)=\left\{\{u_k\}\in\lii;\quad \sum_{k\ge1}k^2|u_k|\le\infty\right\},\qquad \ad a\ e_k = k\ e_k
\end{equation*}

The fundamental commutation rule is $[a;\ad]=\1$

\subsection{The Subharmonic and Pump Fields}

In this paper, $m$ pump fields and $n$ subharmonic fields are considered. These fields are modeled as independent modes of the creator and annihilation operators. The tensor product of $\lii$ denoted by  $\H_p$ and the tensorial product of $\lii$ denoted by $\H_s$ will be considered. The total space is the tensor product between $\H_s$ and $\H_p$, i.e.,

\begin{definition}\label{def:Pump&Sub space}
\[
\begin{array}{c}
\mathrm{Subharmonic\ Space:\ } \H_s :=\bigotimes_{k=1}^n\ell^2(\mathbb{N}), \quad
 \mathrm{Pump\ Space:\ } \H_p :=\bigotimes_{k=1}^m\ell^2(\mathbb{N})\\
\\
 \mathrm{Total\ Space:\ } \H :=\H_s\ot\H_p
\end{array}
\]
\end{definition}

For any finite sequence of natural numbers $\mathsf{s}=(s_i)_{i=1}^n\in\mathbb{N}^n$, it is defined the vector $e(\mathsf{s}):=\ot_{k=1}^ne_{s_i}$, where $e_{s_i}$ indicates the $s_i$-th element of the canonical basis of $\lii$, then $e(\mathsf{s})$ is a generic element of the canonical orthogonal basis of $\H_s$,\\

 Analogously, one can define a generic element of the canonical base of $\H_p$ by $e(\mathsf{p}):=\ot_{k=1}^me_{p_j}$ where $\mathsf{p}=(p_j)_{j=1}^m\in\mathbb{N}^m$\\

Finally, given sequences $\mathsf{s}=(s_i)_{i=1}^n$ and $\mathsf{p}=(p_j)_{j=1}^m$ we define a generic element of the canonical basis of $\H$ by: $e(\mathsf{s},\mathsf{p}):=e(\mathsf{s})\ot e(\mathsf{p})$\\

The following is the formal definition of the subharmonic and pump fields. 

\begin{definition} For $1\le k\le n$ let
\begin{equation*}\begin{array}{c}
\mathrm{Dom}(a_k)=\mathrm{Dom}(\ad_k)= \underbrace{\overbrace{\lii\ot\cdots\ot \mathrm{Dom}(a)}^{k-th}\ot\cdots\ot\lii}_{n-th}\\
\\
a_k:=\underbrace{\overbrace{\1\ot\cdots\ot a}^{k-th}\ot\cdots\1}_{n-th} \qquad\mathrm{and}\qquad
\ad_k:=\underbrace{\overbrace{\1\ot\cdots\ot\ad}^{k-th}\ot\cdots\1}_{n-th}
 \end{array}
\end{equation*}
\end{definition}

The action of the $k$-th subharmonic field over an element of the canonical basis of $\H_s$ is on the $k$-th component of the vector, i.e. formally:
\begin{equation*}\label{eq:field a over base s}
 a_k e(\mathsf{s})=\sqrt{s_k}\ e(\mathsf{s}-\delta_{ik})\quad\mathrm{and}\quad \ad_k e(\mathsf{s})= \sqrt{s_k+1}\ e(\mathsf{s}+\delta_{ik})\quad\mbox{ where  }\mathsf{s}\pm\delta_{ik}=\{s_i\pm\delta_{ik}\}
\end{equation*}

For $1\le k\le m$, an identical construction is possible for pump fields in $\H_p$, thus we get the following
\begin{equation*}\label{eq:field b over base p}
b_ke(\mathsf{p})=\sqrt{p_k}\ e(\mathsf{p}-\delta_{ik})\quad \mathrm{and}\quad \bd_ke(\mathsf{s},\mathsf{p})=\sqrt{p_k+1}\ e(\mathsf{p}+\delta_{ik})
\end{equation*}

\vspace{0.2cm}
\noindent\textbf{REMARK 1} The spaces $\H_s$ and $\H_p$ are included in $\H$ in the natural form, therefore it is possible to identify $a_k$ with $a_k\ot\1$ and $b_k$ with $\1\ot b_k$ on the space $\H$. Therefore, from the equations (\ref{eq:field a over base s}) and \eqref{eq:field b over base p}, it follows that
\begin{equation}\label{eq:field a over base}
a_ke(\mathsf{s},\mathsf{p})=\sqrt{s_k}\ e(\mathsf{s}-\delta_{ik},\mathsf{p})\quad\mbox{and}\quad\ad_ke(\mathsf{s},\mathsf{p})=\sqrt{s_k+1}\ e(\mathsf{s}+\delta_{ik},\mathsf{p})
\end{equation} 
\begin{equation}\label{eq:field b over base}
 b_ke(\mathsf{s},\mathsf{p})=\sqrt{p_k}\ e(\mathsf{s},\mathsf{p}-\delta_{ik})\quad \mathrm{and}\quad \bd_ke(\mathsf{s},\mathsf{p})=\sqrt{p_k+1}\ e(\mathsf{s},\mathsf{p}+\delta_{ik})
\end{equation}

\vspace{0.2cm}
For the $k$-th number operator, of the equation (\ref{eq:field a over base}), it follows that
\begin{equation*}\label{eq:field N over base}
 \ad_ka_k\ e(\mathsf{s})=s_k\ e(\mathsf{s})
\end{equation*}
with $\mathrm{Dom}(\ad_ka_k)= \lii\ot\cdots\ot \mathrm{Dom}(\ad a)\ot\cdots\ot\lii$\\

Through a direct calculation, we obtain the following commutation rules:
\begin{equation}\label{eq:conmmutation-rules}
 [a_k,\ad_k]=[b_k,\bd_k]=\1
\end{equation}
and all the other possible commutations among $a_i,\ad_j, b_k, \bd_l$ are null, for all choices of  $i,j,k$ and $l$.

\section{The Multi-Photon Model}

\subsection{The Hamiltonian Operator}

It is considered that the pump fields arrive with a frequency $w_k^p$ and that subharmonic fields emerge with a frequency $w_k^s$. Due to energy conservation the sum of the frequencies $w_k^p$ and $w_k^s$ must be the same.\\

Due to physical considerations the Hamiltonian term contains free type of energies; the first due to the number of photons pumped $N_p$, the second due to the number of photons emerging $N_s$ and the third to the interaction $I$. The total Hamiltonian is the sum of the three terms:

\begin{equation}\label{eq:Hamiltonian}
H=N_s+N_p+\frac{\i g}{2}I\qquad\mbox{where } g \mbox{ is not a null constant and} 
\end{equation}

\begin{equation}\label{eq:N_s,N_p,I}
N_s=\sum_{i=1}^n w_i^s\ad_i a_i,\quad N_p=\sum_{j=1}^m w_j^p\bd_j b_j,\quad I=\prod_{i=1}^n \ad_i \prod_{j=1}^m b_j-\prod_{i=1}^n a_i\prod_{j=1}^m \bd_j
\end{equation}

\noindent\textbf{REMARK 2}\label{rem:sum w_s=sum w_p}
 It is important to consider the condition: $\ds\sum_{i=1}^n w^s_i=\sum_{j=1}^m w^p_j$.

\subsection{The Channels}

The construction of the mathematical model for the evolution of the system with measurement should include the interaction system-instrument and the loss. These are described by a finite number of channels of the following form:\\
\begin{equation*}\label{eq:Ri-operators}
 R_i^{a}= \alpha_ia_i,\quad R_i^{a^\dag}= \alpha_i^{+}a_i^{\dag},\quad R_j^{b}= \beta_jb_j, \quad\mbox{and}\quad R_j^{b^\dag}= \beta_j^{+}b^{\dag}_b
\end{equation*}

\noindent where $i=1,\dots,n\quad j=1,\dots,m$\ \ and the $\alpha_i$, $\alpha^{+}_i$, $\beta_j$ and $\beta^{+}_j$ are complex numbers. For a physical consideration about the channel you can see [\citen{Bar86LNP} -- \citen{Belavkin-99}], and the final example in [\citen{Castro-Barchielli}].\\

Now, define the operator
\begin{equation}\label{eq:operator R-1}
R:=\sum_{i=1}^n\left(R_i^{a}\right)^*R_i^{a}+\sum_{i=1}^n\left(R_i^{a^\dag}\right)^*R_i^{a^\dag}+ \sum_{j=1}^m\left(R_j^{b}\right)^*R_j^{b} + \sum_{j=1}^m\left(R_j^{b^\dag}\right)^*R_j^{b^\dag}
\end{equation}

>From equation (\ref{eq:operator R-1}) and equation (\ref{eq:conmmutation-rules}), it follows that

\begin{equation*}
 R=\sum_{i=1}^n|\alpha_i|^2\ad a+ \sum_{j=1}^m|\beta_j|^2\bd b+ \left(\sum_{i=1}^n\big|\alpha^{+}_i\big|^2+ \sum_{j=1}^m\big|\beta^{+}_j\big|^2\right)\1
\end{equation*}

\vspace{0.3cm}
\subsection{The Evolution Equation}

In the framework of Quantum Probability theory, the evolution of a  quantum system interacting with a external field is given by a Stochastic Schr\"odinger Equation or Hudson-Parthasaraty Equation\cite{Hudson-Partha, Partha} (H-P equation).\\

The space $\H$, given in Definition \ref{def:Pump&Sub space} is in interaction with a field $\Gamma$, given by the Symmetric Fock space over $\mathrm{L}^2(\mathbb{R}_+,\mathbb{C}^d)$. The elements of the Hudson-Parthasarathy equation can be seen in \emph{Quantum Stochastic Calculus with Unbounded Coefficient} \cite{RCS-LAP} (page 16)

\vspace{0.3cm}
With this definition the Hudson-Parthasarathy equation has the form:

\begin{multline}\label{eq:H-P equation}
  \d U(t)=\displaystyle\Big(\sum_{i\ge1}^n\big(R_i^{a}+R_i^{a^\dag}\big)\d A_i^{\dag}(t)+ \sum_{j\ge1}^m\big(R_j^{b}+R_j^{b^\dag}\big)\d A_i^{\dag}(t)\\
- \sum_{i\ge1}^n\big(R_i^{a}+R_i^{a^\dag}\big)^*\d A_i(t)- \sum_{j\ge1}^m\big(R_j^{b}+R_j^{b^\dag}\big)^*\d  A_j(t)\\ +K\d t\Big)U(t)
\end{multline}
where $K=-iH-\dfrac{1}{2}R$

\vspace{0.3cm}
\section{Existence of the Dynamics}

The conditions for the existence of a solution of the H-P equation with unbounded coefficients were studied in [\citen{Fagnola-Wills}]. In the context of the theory of continuous measurement, these are included In Hypothesis 1 in [\citen{Castro-Barchielli}]. For models such as those described in this paper, these conditions can be reduced to the following theorem.

\begin{theorem}[Theorem 2 in [{\citen{Castro-Barchielli}]}]\label{teo:existence}
Just take $D$ given by the linear span of the basis $\big\{e(\mathsf{s},\mathsf{p});\  \mathsf{s}=\{s_i\}^n_{1},\  \mathsf{p}=\{p_j\}^m_1\subset\mathbb{N}\big\}$ and let $N=N_s+N_p$ where $N_s$ and $N_p$ are the operators given in the equation \eqref{eq:N_s,N_p,I} and consider the following operators
 
$C:=N^{2(n+m)}$, over an appropriate domain and $\Ce:=\dfrac{C}{(1+\epsilon C)^2} $,  $\forall\epsilon>0$\\

If there exist constants $\delta> 0$ and $b_1,\,b_2 \ge 0$ such that $\mr{Dom}(C^{1/2})\subset \mr{Dom}(F)$ and
    
\begin{itemize}
\item for each $\epsilon\in (0,\delta)$, $C_\epsilon^{1/2} {D}\subset D$ and each
    operator $F_{ij}^*C_\epsilon^ {1/2}|_{{D}}$ is bounded.
\item for all $0<\epsilon<\delta$ and $u_0,\ldots,u_d\in \mr{Dom}(F)$, the following
    inequality holds:
    \begin{multline*}
\sum_{i,j\geq0}\left( \la u_i|C_\epsilon F_{ij}u_j\ra +\la F_{ji}u_i|C_\epsilon u_j\ra +
\sum_{k\geq 1}\la F_{ki}u_i|C_\epsilon F_{kj}u_j\ra\right)
\\ {}\leq  \sum_{i\geq 0}\left(b_1
\la u_i| C_\epsilon u_i\ra+ b_2 \norm{u_i}^2\right).
\end{multline*}
 \end{itemize}
 Then, the equation \eqref{eq:H-P equation} admits a unique solution.
\end{theorem}

\begin{definition}\label{def:q}
 Let $e(\mathsf{s},\mathsf{p})$ be a generic element of the canonical basis, we define the functional $q: \H_s\ot\H_p\to\mathbb{C}$ acting over $e(\mathsf{s},\mathsf{p})$ in the following form: 

\begin{equation*}
 q(e(\mathsf{s},\mathsf{p})):=\la e(\mathsf{s},\mathsf{p})\big|Ne(\mathsf{s},\mathsf{p})\ra=\sum_{i=1}^nw_i^ss_i+\sum_{j=1}^mw_i^pp_i
\end{equation*}
\end{definition}

\noindent\textbf{REMARK 3} In order to simplify the notation, we will use only the letter $q$ to indicate the positive number $q(e(\mathsf{s},\mathsf{p}))$\\

>From equations (\ref{eq:field a over base}) and (\ref{eq:field b over base}), for any $ e(\mathsf{s},\mathsf{p})$ element of the basis of $\H$, it is obtained that :
\begin{equation}\label{eq:N over base}
 N e(\mathsf{s},\mathsf{p})=qe(\mathsf{s},\mathsf{p})
\end{equation}

\vspace{0.3cm}
Therefore the number $q$ is the eigenvalue associated to the vector $ e(\mathsf{s},\mathsf{p})$.\\

This allows us to define the operator $f(N)$ in the elements of the basis by:

\begin{equation}\label{eq:fN over base}
f(N) e(\mathsf{s},\mathsf{p}):=f(q) e(\mathsf{s},\mathsf{p}) \qquad\mbox{for any function } f.
\end{equation}

\vspace{0.3cm}
The action of the operator $\Ce$ over the elements of the basis of $\H$, for $\epsilon\ge0$, due to the equations \eqref{eq:N over base} and (\ref{eq:fN over base}), is

\begin{equation}\label{eq:Ce over base}
 \Ce e(\mathsf{s},\mathsf{p})=\frac{C}{(1+\epsilon C)^2}e(\mathsf{s},\mathsf{p})=\frac{q^{2(n+m)}}{(1+\epsilon q^{2(n+m)})^2}e(\mathsf{s},\mathsf{p})
\end{equation}

\begin{proposition}\label{prop:commutation [H,N]}
 The operator $N$, defined in Theorem \ref{teo:existence}, commutes with the Hamiltonian operator $H$, defined in equation (\ref{eq:Hamiltonian}).
\end{proposition}

\noindent\emph{Proof.}
The operator $H=N+I$, then $[H,N]=[N,N]+[I,N]=[I,N]$, therefore, only needs to be proved that $[I,N]=0$

\[
 [I,N]=\left[\prod_{i=1}^n \ad_i \prod_{j=1}^m b_j-\prod_{i=1}^n a_i\prod_{j=1}^m \bd_j\ ,\ N_s+N_p\right]
\]
Analyzing each term, one has:
\[
\begin{array}{rcl}
 \ds\left[\prod_{i=1}^n \ad_i \prod_{j=1}^m b_j\ ,\ N_s\right]&=&\ds\prod_{i=1}^n \ad_i\prod_{j=1}^m b_j \sum_{k=1}^nw_k^s\ad_ka_k-\sum_{k=1}^nw_k^s\ad_ka_k\prod_{i=1}^n \ad_i \prod_{j=1}^m b_j\\
\\
{}&=&\ds\prod_{j=1}^mb_j\sum_{k=1}^nw_k^s\left(\prod_{i\neq1}^n \ad_i\ad_k\ad_ka_k - \ad_ka_k\ad_k\prod_{i\neq1}^n\ad_i\right)\\
\\
{}&=&\ds\prod_{j=1}^mb_j\sum_{k=1}^nw_k^s\prod_{i\neq1}^n\ad_i\ad_k\underbrace{\left(\ad_ka_k-a_k\ad_k\right)}_{-\1}\\
\\
{}&=&\ds-\sum_{k=1}^nw_k^s\left(\prod_{j=1}^mb_j\prod_{i=1}^n \ad_i\right)
\end{array}
\]

Analogously, for the terms: $\ds\left[\prod_{i=1}^n \ad_i \prod_{j=1}^m b_j\ ,\ N_p\right]= \sum_{k=1}^mw_k^p\left(\prod_{j=1}^mb_j\prod_{i=1}^n \ad_i\right)$\\

Therefore, by Remark \ref{rem:sum w_s=sum w_p} it follows that:
 
\[
 \left[\prod_{i=1}^n \ad_i \prod_{j=1}^m b_j\ ,\ N_s\right]+ \left[\prod_{i=1}^n \ad_i\prod_{j=1}^m b_j\ ,\ N_p\right]=0
\]
To proceed in the other two cases, one must proceed in an identical form.\qquad $\square$

\begin{definition}\label{def:Lq}
Now, we define the following auxiliary functionals over elements of the canonical basis of $\H=\H_s\ot\H_p$ by
\[
 \mathscr{L}_0(e(\mathsf{s},\mathsf{p})):=\la e(\mathsf{s},\mathsf{p})\big|C_{\epsilon}e(\mathsf{s},\mathsf{p})\ra\
 \]
\[
\mathscr{L}_k^a(e(\mathsf{s},\mathsf{p})):=\la a_ke(\mathsf{s},\mathsf{p})\big|[C_{\epsilon},a_k]e(\mathsf{s},\mathsf{p})\ra,\quad
 \mathscr{L}_k^{\ad}(e(\mathsf{s},\mathsf{p})):=\la \ad_ke(\mathsf{s},\mathsf{p})\big|[C_{\epsilon},\ad_k]e(\mathsf{s},\mathsf{p})\ra
\]
\[
 \mathscr{L}_k^b(e(\mathsf{s},\mathsf{p})):=\la b_ke(\mathsf{s},\mathsf{p})\big|[C_{\epsilon},b_k]e(\mathsf{s},\mathsf{p})\ra,\quad
 \mathscr{L}_k^{\bd}(e(\mathsf{s},\mathsf{p})):=\la \bd_ke(\mathsf{s},\mathsf{p})\big|[C_{\epsilon},\bd_k]e(\mathsf{s},\mathsf{p})\ra
\]

\end{definition}

\begin{lemma}\label{lemma:L_k in term L_0}

The newly defined functionals can be expressed in terms of the functional $q$ (Definition \ref{def:q}) as follows
\[
 \displaystyle\mathscr{L}_0(e(\mathsf{s},\mathsf{p}))=\dfrac{q^{2(n+m)}}{(1+\epsilon q^{2(n+m)})^2}
\]
\[
 \mathscr{L}_k^a(e(\mathsf{s},\mathsf{p}))=s_k\big(\mathscr{L}_0(q-w_k^s)-\mathscr{L}_0(q)\big),\quad
 \mathscr{L}_k^{\ad}(e(\mathsf{s},\mathsf{p}))=(s_k+1)\big(\mathscr{L}_0(q+w_k^s)-\mathscr{L}_0(q)\big)
\]
\[
 \mathscr{L}_k^b(e(\mathsf{s},\mathsf{p}))=p_k\big(\mathscr{L}_0(q-w_k^p)-\mathscr{L}_0(q)\big),\quad
 \mathscr{L}_k^{\bd}(e(\mathsf{s},\mathsf{p}))=(p_k+1)\big(\mathscr{L}_0(q+w_k^p)-\mathscr{L}_0(q)\big)
\]

\end{lemma}

\noindent\emph{Proof.}
The first equality is immediate. From equation (\ref{eq:Ce over base}) one has 
\[\begin{array}{rcl}
Ne(s_i-\delta_{ik},p_j)&=&\displaystyle\left(\sum_{i\neq k}w_i^ss_i+w_k^s(s_i-1)+\sum_{j=1}^mw_i^pp_i\right)e(s_i-\delta_{ik},p_j)\\
\\
&=&(q-w_k^s)e(s_i-\delta_{ik},p_j).
\end{array}\]

therefore, $\Ce e(s_i-\delta_{ik},p_j)=\mathscr{L}_0(q-w_k^s)e(s_i-\delta_{ik},p_j).$

\[\begin{array}{rcl}
\mathscr{L}_k^a(q)&=&\la a_ke(\mathsf{s},\mathsf{p})\big|[C_{\epsilon},a_k]e(\mathsf{s},\mathsf{p})\ra \\
\\
&=&\la \sqrt{s_k}e(s_i-\delta_{ik},p_j)\big|\Ce a_ke(\mathsf{s},\mathsf{p})-a_k\Ce e(\mathsf{s},\mathsf{p})\ra\\
\\
&=&\la \sqrt{s_k}e(s_i-\delta_{ik},p_j)\big|\Ce\sqrt{s_k}e(s_i-\delta_{ik},p_j)-a_k\mathscr{L}_0(q)e(\mathsf{s},\mathsf{p})\ra\\
\\
&=&\la \sqrt{s_k}e(s_i-\delta_{ik},p_j)\big|\sqrt{s_k}\big(\mathscr{L}_0(q-w_k^s)-\mathscr{L}_0(q)\big)e(s_i-\delta_{ik},p_j)\ra\\
\\
&=&s_k\big(\mathscr{L}_0(q-w_k^s)-\mathscr{L}_0(q)\big)

\end{array}\]
The other cases are analogous. \qquad$\square$\\

The following lemmas, which contains some polynomial inequalities will be useful in proving the main theorems and propositions

\vspace{0.2cm}
\begin{lemma}[Polynomial inequality]\label{lemma:inequality for r&k}
 Let  the constant $r>0$, $M\in\mathbb{N}$ and the polynomial $f(x)=\dfrac{1}{1+\epsilon x^{2M}}$, then for any $x>0$, there exists $\epsilon>0$ such that

\begin{enumerate}
 \item[a)] $f(x-r)\le2^{4M}f(x)$
 \item[b)] $(x-r)^{2M}f(x-r)+x^{2M}f(x)\le2^{2M}x^Mf(x)$
 \item[c)] $s\big|(x-r)^{2M}f(x-r)-x^{2M}f(x)\big|\le2^{2M}x^Mf(x)$, \quad if $s>0$ and $sr<x$
\end{enumerate}

\end{lemma}

\noindent\emph{Proof.}\ \\
 If $x\ge2r$, then $\displaystyle\frac{1}{x-r}\le\displaystyle\frac{2}{x}$  and this implies $\displaystyle\frac{1}{1+\epsilon(x-r)^{2k}} \le\displaystyle\frac{4^k}{1+\epsilon x^{2k}}$ and this implies the result.\\

If $x<2r$ then the inequality is equivalent to $\epsilon\left(x^{2k}-4^k(x-r)^{2k}\right)\le4^k-1$\\

\noindent the maximum of the polynomial $f(x)=x^{2k}-4^k(x-r)^{2k}$ is at 
$x_M=\frac{2^{\frac{2k}{2k-1}}}{2^{\frac{2k}{2k-1}}-1}r$, therefore, it is  $f(x_M)=\displaystyle\frac{2^{\frac{4k^2}{2k-1}}-2^{2k}}{\left(2^{\frac{2k}{2k-1}}-1\right)^{2k}}\ r^{2k}=(2r)^{2k}$\\

Therefore,  $\forall x\in]0,2r[$ one has that $\ds\epsilon\left(x^{2k}-4^k(x-r)^{2k}\right)\le\epsilon (2r)^{2k}$

then just take $\epsilon<\displaystyle\frac{1}{2r^{2k}}$ and the proof of $a)$ is finished.

\vspace{0.2cm}
We provide the inequality b)

\begin{multline*}
  (x-r)^{2M}f(x-r)+x^{2M}f(x)=\left(\dfrac{(x-r)^M}{1+\epsilon(x-r)^{2M}}+ \dfrac{x^M}{1+\epsilon x^{2M}}\right)\\
\\
\le\dfrac{(x-r)^{M}\left(1+\epsilon x^{2M}\right)+x^M\left(1+\epsilon(x-r)^{2M)} \right)}{\left(1+\epsilon(x-r)^{2M}\right)\left(1+\epsilon x^{2M}\right)}\\
 \\
\le\dfrac{\left((x-r)^M+x^M\right)\left(1+\epsilon(x-r)^Mx^M\right)} {\left(1+\epsilon(x-r)^{2M}\right)\left(1+\epsilon x^{2M}\right)}\\
 \\
\le\dfrac{2x^M\left(1+\epsilon x^{2M}\right)}{\left(1+\epsilon(x-r)^{2M} \right)\left(1+\epsilon x^{2M}\right)}=\dfrac{2x^{M}}{1+\epsilon(x-r)^{2M}}\\
\end{multline*}
and part a) implies the result. Analogously, one can show that if $sr<x$, the Binomial Theorem implies that
\[s\left|\frac{(x-r)^M}{1+\epsilon(x-r)^{2M}}-\frac{x^M}{1+\epsilon x^{2M}}\right|
\le s\frac{\left|(x-r)^M-x^M\right|}{1+\epsilon(x-r)^{2M}}
\le\frac{2^{M-1}x^M}{1+\epsilon(x-r)^{2M}}\qquad\square
\]

\begin{lemma}\label{lemma:inequality for L_ab}
 For any $\{s_i\}$ and $\{p_j\}$ one has:

$$\left|\mathscr{L}^a_k(q)\right|\le 32^{n+m}\mathscr{L}_0(q),
\qquad\left|\mathscr{L}^{\ad}_k(q)\right|\le 32^{n+m}\mathscr{L}_0(q)$$

$$\left|\mathscr{L}^b_k(q)\right|\le32^{n+m}\mathscr{L}_0(q),
\qquad\left|\mathscr{L}^{\bd}_k(q)\right|\le 32^{n+m}\mathscr{L}_0(q)$$
\end{lemma}

\noindent\emph{Proof.}\ \\

>From Lemma \ref{lemma:L_k in term L_0}, one has $|\mathscr{L}^a_k(q)|= \big|\mathscr{L}_0(q-w^s_k)-\mathscr{L}_0(q)\big|$ and taking $q=x$, $\displaystyle r=\max_{i,j}\{w^s_i,\ w^p_j \}$ and $\mathscr{L}_0=f$ in Lemma \ref{lemma:inequality for r&k} the result is obtained. The other cases are analogous.\qquad$\square$

\begin{proposition}\label{prop:inequality Re1}
 For any $u\in D$ $\mathrm{Re}\la a_ku\big|[\Ce,a_k]u\ra\le32^{n+m}\la u\big|\Ce u\ra$ and similarly for $\ad_k$, $b_k$ and $\bd_k$.
\end{proposition}

\noindent\emph{Proof.}
 Notice that $\mathrm{Re}\la a_ku\big|[\Ce,a_k]u\ra\le\left|\la a_ku\big|[\Ce,a_k]u\ra\right|$ and we can write the vector $u$ as a sum of the elements of the basis of $\H$. The results are an immediate consequence of Lemma \ref{lemma:inequality for L_ab}\qquad$\square$

\vspace{0.3cm}
\noindent\textbf{REMARK 4}\label{rem:L}
Note that the series $\displaystyle\sum_{(\mathsf{s},\mathsf{p})\in\mathbb{N}^{n+m}}\mathscr{L}_0(q(\mathsf{s},\mathsf{p}))$ converges, so we can define  $\displaystyle\mathscr{L}=\sum_{q}\mathscr{L}_0(q)$
Due to the fact that $\displaystyle\int_{\mathbb{R}^+}\mathscr{L}_0(x)\d x<\infty $, we can define the constant.

\begin{proposition}\label{prop:inequality Re2}
For any choice of $u_0$ and $u$ in $D$ one has
\[
 2\mathrm{Re}\la u\big|[\Ce,a_k]u_0\ra\le32^{n+m} \mathscr{L}\|u\|^2+\|u_0\|^2
\]
where $\mathscr{L}$ is as in REMARK 4.
\end{proposition}

\vspace{0.3cm}
\noindent\emph{Proof.}
Let $\{\gamma^0_l\}$ and $\{\gamma_l\}$ be two families of sequences in $\mathbb{C}$ such that

$ u_0=\sum_{l}\gamma^0_le(\mathsf{s}_l,\mathsf{p}_l)$ and $ u=\sum_{l}\gamma_le(\mathsf{s}_l,\mathsf{p}_l)$,\quad with $\mathsf{s}_l=\{s^l_i\}_1^n$, $\mathsf{p}_l=\{p^l_j\}_1^m$

 \[\begin{array}{rcl}
  2\mathrm{Re}\la u\big|[\Ce,a_k]u_0\ra&=&\ds2\mathrm{Re}\la \sum_{lij}\gamma_le(s_i^l,p_j^l)\Big|[\Ce,a_k]\sum_{l}\gamma^0_le(s_i^l,p_j^l)\ra\\
&=&\ds2\mathrm{Re}\la \sum_{l}\gamma_le(s_i^l,p_j^l)\Big|\sum_{l}\gamma^0_l\frac{\mathscr{L}_k^a(q^l)}{\sqrt{s^l_k}}e(s_i^l-\delta_{ik},p_j^l)\ra\\
&=&\ds2\mathrm{Re}\sum_{l}\gamma_l\overline{\gamma^0_{l+1}}\frac{\mathscr{L}_k^a(q^{l})}{\sqrt{s^l_k+1}}
\end{array}\]

\noindent but note that for any complex numbers $z,w$ one has $2\mathrm{Re}(z\overline{w})\le|z|^2+|w|^2$ and for any even sequences $\{a_i\},\{b_i\}$ of positive numbers one has $\sum_ia_ib_i\le\sum_ia_i\sum_ib_i$. Therefore 
\[
\ds 2\mathrm{Re}\la u\big|[\Ce,a_k]u_0\ra\le 
\sum_{l}\frac{|\mathscr{L}^a_k(q^{l+1})|^2}{s^l_k+1}\sum_{l}\|\gamma_l\|^2+ \sum_{l}\|\overline{\gamma^0_l}\|^2
\]

From Lemma \ref{lemma:inequality for L_ab} one has $\ds\sum_{l}\frac{|\mathscr{L}^a_k(q^{l+1})|^2}{s^l_k+1}\le32^{n+m}\mathscr{L}$,  and the result is obtained.

An analogous result is obtained for $\ad_k, b_k$ and $\bd_k$.\qquad$\square$

\vspace{0.2cm}
Note that if $\displaystyle r=\max_{i,j}\{w_i^s,w_j^p\}$ and $\epsilon<\frac{1}{2r^{2(n+m)}}$, then cleary $C_\epsilon^{\frac{1}{2}}D\subset D$ since all vectors of the canonical basis are eigenvectors of the operator $\Ce^{\frac{1}{2}}$.

\begin{theorem}\label{teo:hipothesis 8a}
The operators $F_{ij}\Ce^{\frac{1}{2}}\big|_D$ are bounded $\forall i,j\ge0$
\end{theorem}
\noindent\emph{Proof.}\\

 Clearly the operators $a_k$, $a^{\dag}_k$, $b_k$ and $b^{\dag}_k$ are relatively bounded with respect to $K=-\i(N_s+N_p+\displaystyle\frac{\i g}{2}I)-\displaystyle\frac{1}{2}R$
 and the component of $K$ ``more'' unbounded is the operator $I$. Therefore if $I\Ce^{\frac{1}{2}}\big|_{D_{\epsilon}}$ is bounded, then $F_{ij}\Ce^{\frac{1}{2}}\big|_{D_{\epsilon}}$ is bounded $\forall i,j\ge0$\\

Let $u\in D$ and ${\gamma}_l$ be a sequence such that $u=\sum_l\gamma_le(\mathsf{s}_l,\mathsf{p}_l)$, with  $\mathsf{s}_l=\{s^l_i\}_1^n$ and $\mathsf{p}_l=\{p^l_j\}_1^m$, then 
\[
 I\Ce^{\frac{1}{2}}u=\sum_l\gamma_l\mathscr{L}_0^{\frac{1}{2}}(q_l)\left(\prod_{i,j}\sqrt{(s^l_i+1)p^l_j}-\prod_{i,j}\sqrt{s^l_i(p^l_j+1)}\right)e' 
\]
 with $e'=\left(e(s^l_i+1,p^l_j-1)-e(s^l_i-1,p^l_j+1)\right)$. Now $s^l_i+1, p^l_j+1\le q_l,\quad \forall i,j$, therefore

\[
 \|I\Ce^{\frac{1}{2}}u\|\le\left\|\sum_l2\gamma_l\frac{q_l^{n+m}}{1+\epsilon q_l^{2(n+m)}} \sqrt{q_l^{n+m}}e'\right\|\le \sum_l2|\gamma_l|\frac{q_l^{\frac{3}{2}(n+m)}}{1+\epsilon q_l^{2(n+m)}}\|e'\|
\]

the sequence $\ds\frac{q_l^{\frac{3}{2}(n+m)}}{1+\epsilon q_l^{2(n+m)}}\to 0$, hence

  $\ds\|I\Ce^{\frac{1}{2}}u\|\le C\sum_l\|\gamma_l\|=C\|u\|$ for some value of the constant $C$.\qquad$\square$

\begin{theorem}\label{theo:existence of the dynamics}
For all $0<\epsilon<\frac{1}{2r^{2k}}$\quad with $\displaystyle r=\max_{i,j}\{w_i^s,w_j^p\}$, there exist constants $b_1,b_2>0$ such that

\[
 2\mathrm{Re}\sum_{i\ge1}\la u_i+\frac{1}{2}R_i^a u_0\Big|[C_{\epsilon},R_i]u_0\ra\le \sum_{i\ge0}\left(b_1\la u_i|C_{\epsilon}u_i\ra+b_2\|u_i\|^2\right)
\]
and analogously for $R_i^{a^{\dag}}$, $R_j^{b}$ and $R_j^{b^{\dag}}$.
\end{theorem}

\noindent\emph{Proof.}

A consequence of Proposition \ref{prop:commutation [H,N]} is that by explicitly computing the left hand side of the inequality (b) of the point (viii) in the Hypothesis 1 in [\citen{Castro-Barchielli}] one has:

\[
 2\mathrm{Re}\sum_{i\ge1}\la u_i+\frac{1}{2}R_i u_0\Big|[C_{\epsilon},R_i]u_0\ra\le \sum_{i\ge0}\left(b_1\la u_i|C_{\epsilon}u_i\ra+b_2\|u_i\|^2\right)
\]

>From equation \eqref{eq:Ri-operators} one has

\[\begin{array}{l}
\displaystyle 2\mathrm{Re}\sum_{k\ge1}\la u_k+\frac{1}{2}R_k^a u_0\Big|[C_{\epsilon},R_k^a]u_0\ra\\
\\
\hspace{1cm}= 2\mathrm{Re}\sum_{k\ge1}\la u_k\Big|[C_{\epsilon},\alpha_ia_k]u_0\ra + \mathrm{Re}\sum_{k\ge1}\la \alpha_ka_k u_0\Big|[C_{\epsilon},\alpha_ka_k]u_0\ra\\
 \end{array}\]

\noindent From Proposition \ref{prop:inequality Re2} 

\[
 2\mathrm{Re}\sum_{k\ge1}\la u_k\Big|[C_{\epsilon},\alpha_ia_k]u_0\ra\le
\sum_{k\ge1}2\mathrm{Re}(\alpha_k)\left(\mathscr{L}\|u_k\|^2+\|u_0\|^2\right)
\]
 
 If $b^{(1)}_2=\mathrm{Max}\big( 2\mathrm{Re}(\alpha_k)\mathscr{L}; 2\mathrm{Re}(\alpha_k)\big)$ then
  $$2\mathrm{Re}\sum_{k\ge1}^n\la u_k\Big|[C_{\epsilon},\alpha_ia_k]u_0\ra\le \sum_{k\ge0}b^{(1)}_2\|u_k\|^2$$

also

\[\begin{array}{rcl}
 \ds\mathrm{Re}\sum_{k\ge1}\la \alpha_ka_k u_0\Big|[C_{\epsilon},\alpha_ka_k]u_0\ra&=&
\ds \sum_{k\ge1}|\alpha_k|^2\mathrm{Re}\la a_ku_0\Big|[C_{\epsilon},a_k]u_0\ra\\
\\
\mbox{By Proposition \ref{prop:inequality Re1}}&\le&\ds \sum_{k\ge1} 32^{n+m}|\alpha_k|^2\la u_0|\Ce u_0\ra\\
\\
&\ds\le& b^{(1)}_1\sum_{k\ge0}\la u_k|\Ce u_k\ra
\end{array}\]

where $\ds b^{(1)}_1:=\sum_{k\ge1} 32^{n+m}|\alpha_k|^2$\\

An analogous procedure is performed for $R^{\ad}, R^b$ and $R^{\bd}$. Then just let us make 
$$b_1=\max_{i=1,\dots,4}\left(b^{(i)}_1\right)\quad \mbox{and}\quad b_2=\max_{i=1,\dots,4}\left(b^{(i)}_2\right)$$ 
and with this the proof is finished.\qquad$\square$

\begin{corollary}\label{coro:existence}
The equation (\ref{eq:H-P equation}) admit a unique solution.
\end{corollary}

\noindent\emph{Proof.}
It is an immediate conclusion of Theorems \ref{teo:existence}, \ref{teo:hipothesis 8a} and \ref{theo:existence of the dynamics}.\qquad$\square$

\subsection{Unitary Solution}

The solution to the equation (\ref{eq:H-P equation}) does not imply, necessarily, that it is unitary. The unitarity of the solution is linked to the Markov Property of the minimal quantum semigroup defined in Theorem 4.5 in [\citen{Fagnola2}]. The unitary property is necessary to have a solution with a physical sense.

\begin{theorem}\label{teo:unitary solution}
 The family of the operators $U(t)$, solution of the equation (\ref{eq:H-P equation}) is a unitary process of the operators.
\end{theorem}

\vspace{0.3cm}
\noindent\emph{Proof.}
 The proof is an application of Theorems 10.2 and 10.3 in [\citen{Fagnola3}], where the operators $A$ and $\ds L_k$ are defined by
$$\ds A=K^*\qquad L_k=-\left(R^a+R^{\ad}+ R^b+R^{\bd}\right)$$
The operators $Q$ and $Z$ are defined by.

$$Q=k_1(N_s+N_p) + k_2\1 \quad\mbox{and}\quad Z=R\qquad\mbox{where}$$
$$k_1=\max_{k}\{|\alpha_k|^2,|\alpha^+_k|^2, |\beta_k|^2,|\beta^+_k|^2\} \quad\mbox{and}\quad  k_2=\sum_{k}|\alpha_k^+|^2+|\beta_k^+|^2$$
The proof continues as in [\citen{Castro-Barchielli}, Proposition 6.1]
\begin{flushright}$\square$\end{flushright}

\end{document}